\documentclass[12pt,thmsa]{article}
\usepackage{sw20lart}

\input{tcilatex}
\begin{document}

\title{Illumination by Taylor Polynomials}
\author{Alan Horwitz}
\date{(Received 22 November 1999)}
\maketitle

\begin{abstract}
Let $f(x)$ be a differentiable function on the real line $R$, and let $P$ be
a point not on the graph of $f(x)$. Define the illumination index of $P$ to
be the number of distinct tangents to the graph of $f$ which pass thru $P$.
We prove that if $f^{\prime \prime }$ is continuous and nonnegative on $R$, $%
f^{\prime \prime }\geq m>0$ outside a closed interval of $R$, and $f^{\prime
\prime }$ has finitely many zeros on $R$, then any point $P$ \textit{below }%
the graph of $f$ has illumination index $2$. This result fails in general if 
$f^{\prime \prime }$ is not bounded away from $0$ on $R$. Also, if $%
f^{\prime \prime }$ has finitely many zeros and $f^{\prime \prime }$ is not
nonnnegative on $R$, then some point below the graph has illumination index
not equal to $2$. Finally, we generalize our results to illumination by odd
order Taylor polynomials.

\end{abstract}

\section{\protect\smallskip Introduction}

The central problem in differential calculus is to find the tangent line to
a given curve $y=f(x)$ at a given point $(c,f(c))$ \textit{on the graph of }$%
f$. A somewhat more complicated problem is: Given a point $P=(s,t)$ \textit{%
not on the graph of }$f$, find all values of $c$ so that the tangent line to
the graph of $f$ at $(c,f(c))$ passes through $P.$ If such a $c$ exists, we
say that the point $(c,f(c))$ \textit{illuminates} $P$. A typical example
is: find all tangents to $y=x^{2}$ which pass through the point $(2,3)$. In
this case each of the points $(1,1)$ and $(3,9)$ would illuminate $P$. Of
course, it is certainly possible that \textit{no} tangent line at all passes
through the given point $(s,t)-$e.g. if $y=x^{2}$ and $P=(1,3).$ A simple,
but interesting exercise is: Let $P$ be any point below the graph of $%
y=x^{2} $. Prove that there are exactly two tangents to the graph which pass
through $P$. In considering this type of problem, the following question
naturally arises: given $f(x)$, for which points $P=(s,t)$ is there a
tangent line to $f$ which passes through $P$ ? Also, how many tangents pass
through $P$ ? \smallskip The questions above lead to some potentially
interesting ideas for research. For example, suppose that $f$ is convex on $%
R $, and let $P$ be any point below the graph of $y=f(x)$. Are there always
exactly two tangents to the graph which pass through $P$ ? What if one
assumes that $f^{\prime \prime }(x)>0$ on $R$ ? We give the answers in
Section 1(see, in particular, Theorem \ref{T1}).

In Section 2, we prove a converse result to Theorem \ref{T1}. It is also
natural to try to extend our results to illumination by \textit{higher }%
order Taylor polynomials. In Section 3 we prove results similar to Theorem 
\ref{T1} for illumination by \textit{odd }order Taylor polynomials. Most of
the proofs extend verbatim, but some results from \cite{2} are needed.

\section{Illumination by Tangent Lines}

\begin{definition}
Let $f(x)$ be a differentiable function on the real line, and let $P$ be any
point not on the graph of $f$. We say that the \textit{illumination index}
of $P$ is $k$ if there are $k$ distinct tangents to the graph of $f$ which
pass through $P$. We include the possibility that $k=\infty .$
\end{definition}

\begin{remark}
We say that a tangent line $T$ is \textit{multiple} if $T$ is tangent to the
graph of $f$ at more than one point. If only one tangent line $T$ passes
through $P$, but $T$ is a multiple tangent, we still define the illumination
index of $P$ to be one. One could, of course, define an illumination index
which takes into account the number of points of tangency of each tangent
line.
\end{remark}

As noted earlier, any point below the graph of $y=x^{2}$ has illumination
index $2$. We now generalize this to convex $C^{2}$ functions in general,
with the added condition that $f^{\prime \prime }$ is bounded below by a
positive number outside some closed interval(see Theorem \ref{T1} below).
First we prove a couple lemmas.

\begin{lemma}
\label{L1}Let $f(x)\in C^{2}(R)$, and suppose that there exists $T>0$ such
that $f^{\prime \prime }(x)\geq m>0$ on $\left| x\right| >T$. Let $T_{c}(x)$
be the tangent line to $f$ at $(c,f(c))$. Then for any fixed $%
s,\lim\limits_{\left| c\right| \rightarrow \infty }T_{c}(s)=-\infty $.
\end{lemma}

\smallskip 
%TCIMACRO{\TeXButton{Proof}{\proof} }
%BeginExpansion
\proof%
%EndExpansion
For fixed $s$, let $g(c)=T_{c}(s),$ which implies that $g^{\prime
}(c)=(s-c)f^{\prime \prime }(c)$. Let $U=\max (s,T),$ $u=\min (s,-T)$. It
follows that $g^{\prime }(c)$ is $\QDATOPD\{ \} {\leq \;0,\;c>U}{\geq
\;0,\;c<u}\Rightarrow g(c)$ is $\QDATOPD\{ \} {\text{decreasing on }%
(U,\infty )}{\text{increasing on }(-\infty ,u)}.$ Also, since $f^{\prime
\prime }(x)\geq m>0$ on $\left| x\right| >T,$ 
\begin{equation}
\QATOPD\{ \} {\lim\limits_{c\rightarrow \infty }g^{\prime }(c)=-\infty
}{\lim\limits_{c\rightarrow -\infty }g^{\prime }(c)=\infty }  \label{e1}
\end{equation}

Partition $[U,\infty )$ into infinitely many subintervals, $[c_{k-1},c_{k}],$
of constant width $h>0$. By (\ref{e1}), given $M>0$, there exists $C>0$ such
that $g^{\prime }(c)\leq -M$ for $c\geq C$. Now $g(c_{k})=g(c_{k-1})+\dint%
\limits_{c_{k-1}}^{c_{k}}g^{\prime }(t)dt\leq g(U)-Mh$ if $c_{k-1}\geq C$.
Since this inequality holds for any $M>0,$ $g(c_{k})\rightarrow -\infty $.
Also, since the inequality holds for any increasing sequence $%
\{c_{k}\}\rightarrow \infty ,$ with $c_{k}-c_{k-1}$ constant, $%
\lim\limits_{c\rightarrow \infty }g(c)=-\infty $. A similar argument shows
that $\lim\limits_{c\rightarrow -\infty }g(c)=-\infty $. 
%TCIMACRO{\TeXButton{End Proof}{\endproof}}
%BeginExpansion
\endproof%
%EndExpansion

\begin{remark}
Lemma \ref{L1} is a little easier to prove under the stronger assumption
that $f^{\prime \prime }(x)$ is positive and bounded away from $0$ on the
real line. One can then just examine the error $E_{c}(x)=f(x)-T_{c}(x)$ and
use Taylor's Remainder formula.
\end{remark}

\begin{lemma}
\label{L2}Suppose that $f^{\prime \prime }(x)$ is continuous, nonnegative,
and has \textit{finitely many }zeros in $R$. Then at most two distinct
tangent lines to $f$ can pass through any given point $P$ in the plane.
\end{lemma}

%TCIMACRO{\TeXButton{Proof}{\proof}}
%BeginExpansion
\proof%
%EndExpansion
Suppose that three distinct tangents, $T_{1},T_{2},T_{3}$ pass through $P,$
and suppose that the $T_{i}$ are tangent at $(x_{i},f(x_{i})),i=1,2,3.$
Assume, without loss of generality, that $x_{1}<x_{2}<x_{3}.$ Since $f$ is
convex on any open interval, each pair of tangents has a unique point of
intersection. Let $I_{1}$ equal\ the intersection point of $T_{1}$ and $%
T_{2},$and let $I_{2}$ equal\ the intersection point of $T_{2}$ and $T_{3}.$
Since all three tangents pass through $P$, $I_{1}=I_{2}=P$. If $%
I_{1}=(s_{1},t_{1})$ and $I_{2}=(s_{2},t_{2}),$ then, again, since $f$ is
convex on any open interval, $x_{1}<s_{1}<x_{2}$ and $x_{2}<s_{2}<x_{3},$
which implies that $s_{1}<s_{2},$which contradicts the fact that $%
I_{1}=I_{2}.$

\begin{theorem}
\label{T1}Suppose that $f^{\prime \prime }(x)$ is continuous, nonnegative,
and has \textit{finitely many }zeros in $R$. Assume also that there exists $%
T>0$ such that $f^{\prime \prime }(x)\geq m>0$ on $\left| x\right| >T$. Let $%
P=(s,t)$ with $t<f(s{})$. Then there are exactly two distinct tangent lines
to the graph of $f$ which pass through $P$.
\end{theorem}

%TCIMACRO{\TeXButton{Proof}{\proof}}
%BeginExpansion
\proof%
%EndExpansion
Since $t<f(s{}),$ for $c\;$sufficiently close to $s,T_{c}(s)$ $%
=f(c)+f^{\prime }(c)(s-c)>t.$ By Lemma \ref{L1}, $\lim\limits_{\left|
c\right| \rightarrow \infty }T_{c}(s)=-\infty .$ Hence, for $\left| c\right| 
$ sufficiently large, $T_{c}(s)<t.$ By the Intermediate Value Theorem, $%
T_{c}(s)=t$ for at least two values of $c.$ Note also that for a convex
function, $c_{1}\neq c_{2}\Rightarrow T_{c_{1}}\neq T_{c_{2}}$. Hence the
illumination index of $P$ is at least two. By Lemma \ref{L2}, the
illumination index of $P$ is at most two. That proves the theorem. 
%TCIMACRO{\TeXButton{End Proof}{\endproof}}
%BeginExpansion
\endproof%
%EndExpansion

\medskip The following example shows that Theorem \ref{T1} does \textit{not }%
hold in general for functions which only satisfy $f^{\prime \prime }(x)>0$
on $R$.

\begin{example}
Let $f(x)=\dint\limits_{0}^{x}\left( \dint\limits_{0}^{t}e^{-u^{2}}du\right)
dt=\dfrac{1}{2}\func{erf}\left( x\right) \sqrt{\pi }x+\frac{1}{2}e^{-x^{2}}-%
\dfrac{1}{2},$where $\func{erf}(x)=\dint\limits_{0}^{x}e^{-t^{2}}dt.$ Since $%
f^{\prime \prime }(x)=e^{-x^{2}}$, $\lim\limits_{x\rightarrow \pm \infty
}f^{\prime \prime }(x)=0.$ We now show that \textit{no} tangent to $f$
passes through the point $(0,t)$ when $t<-\dfrac{1}{2}<f(0)=0.$ If the
tangent line $T_{c}$ to $f$ passes through $(s,t),$ then $f(c)+f^{\prime
}(c)(s-c)=t.$ So consider the function $h(x)=f(x)+f^{\prime }(x)(s-x)-t=%
\dfrac{1}{2}e^{-x^{2}}-\dfrac{1}{2}+\dfrac{1}{2}\func{erf}\left( x\right) 
\sqrt{\pi }s-t.$ If $s=0,$ then $h(x)=\dfrac{1}{2}e^{-x^{2}}-\dfrac{1}{2}%
-t\Rightarrow $ $h^{\prime }(x)=-xe^{-x^{2}},$ which implies that $h(x)$ is
increasing for $x<0$ and decreasing for $x>0.$ Since $\lim\limits_{x%
\rightarrow \pm \infty }\left( \dfrac{1}{2}e^{-x^{2}}-\dfrac{1}{2}-t\right)
=\allowbreak -\dfrac{1}{2}-t,$ if $t<-\dfrac{1}{2}$ then $h$ is always
positive and thus has no real zeroes.
\end{example}

Our definition of the illumination index $k$ includes the possibility that $%
k=\infty .$ Of course, for polynomials the illumination index is always
finite(indeed, it's bounded above by the degree of the polynomial). The
following example shows that there are entire functions, however, where 
\textit{almost} \textit{every} point not on the graph has infinite
illumination index.

\begin{example}
Let $f(x)=\sin x$, and let $P=(s,t)$ be any point not on the graph of $f$,
with $t\neq \pm 1$. The tangent line at $(c,f(c))$ passes through $P$ if and
only if $f(c)+f^{\prime }(c)(s-c)=t$, that is, when $g(c)=\sin c+(s-c)\cos
c-t=0$. For $n$ sufficiently large and even, $g(n\pi )=(-1)^{n}(s-n\pi )-t<0$%
, while for $n$ sufficiently large and odd, $g(n\pi )>0$. Hence $g$ has
infinitely many zeroes $c_{1},c_{2},c_{3},...$ Note that since $t\neq \pm 1$%
, none of the zeroes is an odd multiple of $\dfrac{\pi }{2}$, and hence none
of the tangents at $(c_{j},\sin c_{j})$ is horizontal. Now each of these
tangents passes through $P$, but they may not all be \textit{distinct}.
However, since a nonhorizontal line can only be tangent to $y=\sin x$ at
finitely many points, it is clear that infintely many distinct tangents pass
through $P$, and thus $P$ has infinite illumination index.
\end{example}

\begin{remark}
Given $f$, one may define, for each nonnegative integer $k$, the set $D_{k}$%
, equal to the set of points in the plane with illumination index $k$. The $%
D_{k}$ form a partition of $R^{2}-G$, where $G$ is the graph of $f.$ For
example, if $f(x)=x^{3}$, it is not hard to show that $D_{3}=%
\{(s,t):s>0,0<t<s^{3}\}\cup \{(s,t):s<0,s^{3}<t<0\},D_{2}=\{(s,t):s\neq
0,t=0\},D_{1}=G-(D_{2}\cup D_{3})$, and $D_{k}=\emptyset $ for $k=0$ or $k>3$%
.
\end{remark}

\section{A Converse Result}

Suppose that $f(x)$ is \textit{not} convex on $R$. Is it possible for every
point below the graph of $f$ to have illumination index $2$ ? The answer is
no, and thus we have the following partial converse of Theorem \ref{T1}.

\begin{theorem}
\label{T2}Let $f\in C^{3}(R)$ and suppose that $f^{\prime \prime }(x)$ has
finitely many zeroes in $R$. If $f^{\prime \prime }(x)$ is \textit{not}
nonnegative on $R,$ then there is a point $P$ below the graph of $f$ with
illumination index \textit{not} equal to $2$.
\end{theorem}

%TCIMACRO{\TeXButton{Proof}{\proof}}
%BeginExpansion
\proof%
%EndExpansion
. If $f^{\prime \prime }(x)\leq 0$ on $R$, then clearly any point $P$ below
the graph of $f$ has illumination index $0$. Hence we may assume that there
are real numbers $s$ and $u$ such that $f^{\prime \prime }(s)>0,f^{\prime
\prime }(u)=0$, and $f^{\prime \prime }(x)$ changes sign at $x=u$, with $%
f^{\prime \prime }\geq 0$ between $s$ and $u.$ We consider the case $s<u$,
the other case being similar. Let $P=(s,t)$, with $t$ to be chosen shortly.
Now the tangent line at $(c,f(c))$ passes through $P$ if and only if $%
f(c)-cf^{\prime }(c)+f^{\prime }(c)s=t$, which holds if and only if $h(c)=0$%
, where 
\[
h(c)=f(c)-cf^{\prime }(c)+f^{\prime }(c)s 
\]

$h^{\prime }(c)=(s-c)f^{\prime \prime }(c)$ and $h^{\prime \prime
}(c)=(s-c)f^{\prime \prime \prime }(c)-f^{\prime \prime }(c)$. Note that $%
h(s)=f(s)$, $h^{\prime }(s)=0,$ and $h^{\prime \prime }(s)=-f^{\prime \prime
}(s)<0$, so that $h(s)$ is a local maximum of $h(c)$. Since $h^{\prime
}(c)\leq 0$ on $(s,u)$ and $h^{\prime }(c)\geq 0$ on $(u,u+\epsilon )$, $%
h(u) $ is a local minimum of $h(c)$. Note that $h(u)<h(s).$ Let $T$ be the
line $y=h(u)$, the tangent to $h$ at $(u,h(u))$.

\textbf{Case 1:} $T$ only intersects the graph of $h$ at $(u,h(u)).$

Then let $t=h(u).$

\textbf{Case 2:} $T$ intersects the graph of $h$ at some point $Q\neq
(u,h(u)).$ If a $Q$ exists such that $h-T$ changes sign at $Q$, then $%
y=h(u)+\epsilon $ intersects the graph of $h$ in at least three points for
some $\epsilon >0$. If no such $Q$ exists, then $h$ must have another local
minimum at $Q$. Then $y=h(u)+\epsilon $ intersects the graph of $h$ in at
least four points for some $\epsilon >0$. In either case, let $%
t=h(u)+\epsilon $, with $\epsilon $ chosen sufficiently small so that $%
h(u)+\epsilon <h(s).$ Since the zeroes of $h$ correspond to values of $c$
such that the tangent line at $(c,f(c))$ passes through $P$, for case two
there are at least three such values of $c.$ However, it is possible that
some of the corresponding tangents could be multiple. It was shown in \cite
{3}, however, that $f$ can have only \textit{finitely many} multiple tangent
lines in any bounded interval. Also, since each tangent is tangent at only
finitely many points, we can also choose $\epsilon $ sufficiently small so
that none of the tangents corresponding to the zeroes of $h$ is multiple.
Thus at least three \textit{distinct }tangents pass through $P$.

In each case covered, $P$ lies below the graph of $f$ since $h(s)=f(s)$.
Hence the illumination index of $P$ is either one or greater than or equal
to three, and thus cannot equal two. 
%TCIMACRO{\TeXButton{End Proof}{\endproof}}
%BeginExpansion
\endproof%
%EndExpansion

\section{Illumination by Higher Order Taylor Polynomials}

The results of the previous section can be extended to illumination by
Taylor polynomials of order $r$, $r$ odd. In certain ways, the odd order
Taylor polynomials $P_{c}(x)$ behave like tangent lines. Suppose that $f\in
C^{r+1}(-\infty ,\infty ),$ and let $P_{c}(x)$ denote the Taylor polynomial
to $f$ of order $r$ at $x=c$. In \cite{1} it was proved that if $%
f^{(r+1)}(x)\neq 0$ on $[a,b]$, then there is a unique $u,\;a<u<b,$ such
that $P_{a}(u)=P_{b}(u)$. This defines a mean $m(a,b)\equiv u$. We shall
prove a slightly stronger version of this result. The method of proof is
very similar to that used in \cite{2}$,$ where further results and
generalizations of the means $m(a,b)$ were proved.

For the rest of this section we assume that $r$ is an \textit{odd} positive
integer.

Let $E_{c}(x)=f(x)-P_{c}(x)$. By the integral form of the remainder, we have

\begin{equation}
E_{c}(x)=\frac{1}{r!}\dint\limits_{c}^{x}f^{(r+1)}(t)(x-t)^{r}dt  \label{eq3}
\end{equation}

\begin{lemma}
\label{lm1}Suppose that $f^{(r+1)}(x)$ is continuous, nonnegative, and has 
\textit{finitely many }zeros in $[a,b]$. Then $P_{b}-P_{a}$ has precisely
one real zero $c,\;a<c<b.$
\end{lemma}

%TCIMACRO{\TeXButton{Proof}{\proof}}
%BeginExpansion
\proof%
%EndExpansion
By (\ref{eq3}),

\[
E_{a}(x)=\frac{1}{r!}\dint\limits_{a}^{x}f^{(r+1)}(t)(x-t)^{r}dt 
\]

\[
E_{b}(x)=\frac{1}{r!}\dint\limits_{x}^{b}f^{(r+1)}(t)(t-x)^{r}dt 
\]

$E_{a}^{\prime }(x)=\dfrac{1}{(r-1)!}\dint%
\limits_{a}^{x}f^{(r+1)}(t)(x-t)^{r-1}dt\Rightarrow E_{a}^{\prime }(x)<0$
for $x<a$ and $E_{a}^{\prime }(x)>0$ for $x>a.$ Hence $E_{a}(x)$ is strictly
increasing on $(a,b)$. Similarly, $E_{b}(x)$ is strictly decreasing on $%
(a,b) $. Since $E_{a}(a)=0$ and $E_{b}(b)=0,$ there is a unique $c,\;a<c<b,$
such that $E_{b}(c)-E_{a}(c)=0.$ This implies that $P_{b}(c)-P_{a}(c)=0$. Now

$(E_{b}-E_{a})(x)=-\dint\limits_{a}^{b}f^{(r+1)}(t)(x-t)^{r}dt\Rightarrow
(E_{b}-E_{a})^{\prime
}(x)=-r\dint\limits_{a}^{b}f^{(r+1)}(t)(x-t)^{r-1}dt\leq 0$ for $x\in R$.
Since $f^{(r+1)}$ has finitely many zeros, this implies that $E_{b}-E_{a}$
is strictly decreasing on $R$. Hence $E_{b}-E_{a}\;$has precisely one real
zero, which implies that $P_{b}-P_{a}$ has precisely one real zero $%
c,\;a<c<b.$

\begin{lemma}
\label{lm2}Suppose that $f^{(r+1)}(x)$ is continuous, nonnegative, and has 
\textit{finitely many }zeros in $[a,b]$. Let $P$ be any point in the $xy\;$%
plane. Then at most two distinct Taylor polynomials of order $r$ at $%
x=c,\;a\leq c\leq b,$ can pass through $P$.
\end{lemma}

%TCIMACRO{\TeXButton{Proof}{\proof}}
%BeginExpansion
\proof%
%EndExpansion
Suppose that three distinct Taylor polynomials of order $r$, $%
P_{c_{1}},P_{c_{2}},$ and $P_{c_{3}},$ pass through $P=(s,t)$. Then $%
(P_{c_{2}}-P_{c_{1}})(s)=0$ and $(P_{c_{3}}-P_{c_{2}})(s)=0$. Without loss
of generality, assume that $a\leq c_{1}<c_{2}<c_{3}\leq b$. By Lemma \ref
{lm1}, $c_{1}<s<c_{2}$ and $c_{2}<s<c_{3}$, which is a contradiction. Hence
at most two distinct Taylor polynomials of order $r$ can pass through $P$.

\begin{lemma}
\label{lm3}Let $f(x)\in C^{r+1}(R).$ Suppose that there exists $T>0$ such
that $f^{(r+1)}(x)\geq m>0$ on $\left| x\right| >T$. Then for any fixed $%
s,\lim\limits_{\left| c\right| \rightarrow \infty }P_{c}(s)=-\infty $.
\end{lemma}

\smallskip 
%TCIMACRO{\TeXButton{Proof}{\proof} }
%BeginExpansion
\proof%
%EndExpansion
The proof is almost identical to that of Lemma \ref{L1}, and we omit it. 
%TCIMACRO{\TeXButton{End Proof}{\endproof}}
%BeginExpansion
\endproof%
%EndExpansion

\begin{theorem}
\label{Tr2}Suppose that $f^{(r+1)}(x)$ is continuous, nonnegative, and has 
\textit{finitely many }zeros in $R$. In addition, assume that there exists $%
T>0$ such that $f^{(r+1)}(x)\geq m>0$ on $\left| x\right| >T$. Let $P=(s,t)$
with $t<f(s{})$. Then there are exactly two distinct Taylor polynomials of
order $r$ to the graph of $f$ which pass through $P$.
\end{theorem}

%TCIMACRO{\TeXButton{Proof}{\proof}}
%BeginExpansion
\proof%
%EndExpansion
Since $t<f(s{}),$ for $c\;$sufficiently close to $s,$

$P_{c}(s)$ $=f(s)+\dsum\limits_{k=1}^{r}\dfrac{f^{(k)}(c)}{k!}(s-c)^{k}>t.$
By Lemma \ref{lm3}, $\lim\limits_{\left| c\right| \rightarrow \infty
}P_{c}(s)=-\infty $, and hence, for $\left| c\right| $ sufficiently large, $%
P_{c}(s)<t.$ By the Intermediate Value Theorem, $P_{c}(s)=t$ for at least
two values of $c.$ Also, it is not hard to show that if $f^{(r+1)}(x)>0$ on $%
R$, then $c_{1}\neq c_{2}\Rightarrow P_{c_{1}}\neq P_{c_{2}}$. Hence the
illumination index of $P$ is at least two. By Lemma \ref{lm2}, it is \textit{%
at most} two. That proves the theorem. 
%TCIMACRO{\TeXButton{End Proof}{\endproof}}
%BeginExpansion
\endproof%
%EndExpansion

\begin{example}
Let $f(x)=e^{x}+x^{4},P=(0,0),r=3$. Then Theorem \ref{Tr2} applies, and the
illumination index of $P$ equals $2$. We now verify this by estimating the
actual values of $c$. Let $P_{c}(x)$ be the third order Taylor polynomial to 
$f$ at $(c,f(c))$. Then $P_{c}(0)=\allowbreak e^{c}-e^{c}c+\dfrac{1}{2}%
e^{c}c^{2}-\dfrac{1}{6}e^{c}c^{3}-c^{4}=0$ has solutions $c_{1}\approx $ $%
-.9953$ and $c_{2}\approx .\,9782$. Note that if $f(x)=e^{x}$ instead, then
the illumination index of $P$ equals $1$. This does not contradict Theorem 
\ref{Tr2} since $f^{(iv)}(x)\rightarrow 0$ as $x\rightarrow -\infty $.
\end{example}

\end{document}